\documentstyle[amscd]{amsart}
\newtheorem{lemma}{Lemma}[section]
\newtheorem{proposition}[lemma]{Proposition}

\newtheorem{theorem}[lemma]{Theorem}
\newtheorem{corollary}[lemma]{Corollary}

\numberwithin{equation}{section}

\def\NN{{\mathbb N}}
\def\PP{{\mathbb P}}

\def\ZZ{{\mathbb Z}}

\def\Mol{{\bar M}}

\def\Ext{\operatorname {Ext}}

\def\Hom{\operatorname {Hom}}

\def\1I{\operatorname {I}}
\def\2I{\operatorname {II}}

\def\th{\operatorname {th}}    

\def\dim{\operatorname{dim}}

\def\Ext{\operatorname{Ext}}
\def\Fdim{{\sf Fdim}}

\def\Fract{\operatorname{Fract}}

\def\Gr{{\sf Gr}}
\def\GrMod{{\sf GrMod}}

\def\Hom{\operatorname{Hom}}

\def\id{\operatorname{id}}

\def\liminj{\varinjlim}

\def\Mod{{\sf Mod}}

\def\Proj{\operatorname{Proj}}
\def\Projnc{\operatorname{Proj}_{nc}}

\def\Spec{\operatorname{Spec}}

\def\Tails{{\sf Tails}}

\def\l{\leftarrow}

\def\d{\downarrow}

\def\a{\alpha}
\def\b{\beta}

\def\d{\delta}

\def\l{\lambda}

\def\s{\sigma}

\def\fm{{\mathfrak m}}

\def\sA{{\sf A}}
\def\sB{{\sf B}}

\def\sS{{\sf S}}
\def\sT{{\sf T}}

\def\cL{{\cal L}}
\def\cM{{\cal M}}
\def\cN{{\cal N}}
\def\cO{{\cal O}}

\def\Qcoh{{\sf Qcoh}}

\newdimen\uboxsep \uboxsep=1ex
\def\uboxn#1{\vtop to 0pt{\hrule height 0pt depth 0pt\vskip\uboxsep
\hbox to 0pt{\hss #1\hss}\vss}}

\def\uboxs#1{\vbox to 0pt{\vss\hbox to 0pt{\hss #1\hss}
\vskip\uboxsep\hrule height 0pt depth 0pt}}

\def\hoek{\hbox{\vtop{\hbox{\vrule\phantom{xx}\vrule}\hrule}\kern -0.4pt}}

\def\dirlim{\mathop{\vtop{\baselineskip -100pt\lineskip -1pt\lineskiplimit 0pt
\setbox0\hbox{lim}\copy0\hbox to \wd0{\rightarrowfill}}}\limits}
\def\invlim{\mathop{\vtop{\baselineskip -100pt\lineskip -1pt\lineskiplimit 0pt
\setbox0\hbox{lim}\copy0\hbox to \wd0{\leftarrowfill}}}\limits}

\def\I11{{1 \kern -0.8pt \! \mbox{l}}}
\def\mumu{{\mu\kern-4.2pt\mu}}

\def\boxtimes{\setbox0\hbox{$\Box$}\copy0\kern-\wd0\hbox{$\times$}}

\def\barotimes{\operatorname{\bar{\otimes}}}
\def\uuHom{\operatorname{\underline{\underline{Hom}}}}
\newtheorem{example}[lemma]{Example}

\begin{document}
\pagenumbering{roman}

\title{Maps between non-commutative spaces}
\author{S. Paul Smith} 
\address{Department of Mathematics, Box 354350, Univ.
Washington, Seattle, WA 98195, USA}
\email{smith@@math.washington.edu}

\thanks{The author was supported by NSF grant DMS-0070560 }

\begin{abstract}
Let $J$ be a graded ideal in a not necessarily commutative
graded $k$-algebra $A=A_0 \oplus A_1 \oplus \cdots$ in which
$\dim_k A_i < \infty$ for all $i$.
We show that the map $A \to A/J$ induces a closed immersion 
$i:\Projnc A/J \to \Projnc A$ between the non-commutative
projective spaces with homogeneous coordinate rings $A$ and $A/J$. 
We also examine two other kinds of maps between non-commutative spaces.
First, a homomorphism $\phi:A \to B$ between not necessarily
commutative $\NN$-graded rings, induces an affine map 
$ \Projnc B \supset U  \to \Projnc A $ 
from a non-empty open subspace $U \subset \Projnc B$.
Second, if $A$ is a right noetherian connected graded algebra (not
necessarily generated in degree one), and $A^{(n)}$ is a Veronese
subalgebra of $A$, there is a map $\Projnc A \to \Projnc
A^{(n)}$; we identify open subspaces on which this map is an 
isomorphism. Applying these general results when $A$ is (a quotient
of) a weighted polynomial ring produces
a non-commutative resolution of (a closed subscheme of) a
weighted projective space.
\end{abstract}

\maketitle

\pagenumbering{arabic}

\section{Introduction}

Following Rosenberg \cite{R} and Van den Bergh 
\cite{vdB}, a non-commutative space $X$ is a
Grothendieck category $\Mod X$. A {\sf map} $g:Y \to X$ between two
spaces is an adjoint pair of functors $(g^*,g_*)$ with $g_*:\Mod Y
\to \Mod X$ and $g^*$ left adjoint to $g_*$. 
The map $g$ is {\sf affine} \cite[page 278]{R} if $g_*$ is 
faithful and has a right adjoint.
For example, a ring homomorphism $\varphi:R \to S$ induces 
an affine map $g:Y \to X$ between the affine spaces defined by
$\Mod Y:=\Mod S$ and $\Mod X:=\Mod R$.

We say that $g:Y \to X$ is a {\sf closed immersion} if it is affine
and the essential image of $\Mod Y$ in $\Mod X$ under $g_*$ is
closed under submodules and quotients. 

This paper concerns maps between non-commutative projective spaces 
of the form $\Projnc A$. We recall the definition.

Let $k$ be a field. An $\NN$-graded $k$-algebra $A$ is {\sf locally
finite} if $\dim_k A_i < \infty$ for all $i$. The 
non-commutative projective space $X$ with homogeneous 
coordinate ring $A$ is defined by
$$
\Mod X:= \GrMod A/\Fdim A
$$
(see Section \ref{sect.defns}), and
$$
\Projnc A := (\Mod X, \cO_X),
$$
where $\cO_X$ is the image of $A$ in $\Mod X$.
Thus $\Projnc A$ is an enriched quasi-scheme in the languauge of
\cite{vdB}. Let $Y$ be another non-commutative projective space
with homogeneous coordinate ring $B$. A map $f:\Projnc  B \to \Projnc
A$ is a map $f:Y \to X$ such that $f^* \cO_X \cong \cO_Y$.

When $A$ is a commutative $\NN$-graded $k$-algebra we write 
$\Proj A$ for the usual projective scheme. We will always view a
quasi-separated, quasi-compact scheme $X$ as a non-commutative space by
associating to it the enriched space $(\Qcoh X,\cO_X)$. The rule $X
\mapsto (\Qcoh X,\cO_X)$ is a faithful functor.

The main results in this paper are 
Theorems \ref{thm.closed.proj}, \ref{thm.ratl.map}, \ref{thm.Ver},
and Proposition \ref{prop.open.isom}.

Theorem \ref{thm.closed.proj} shows that a surjective homomorphism
$A \to A/J$ of graded rings induces a closed immersion 
$i:\Projnc A/J \to \Projnc A$. The functors $i^*$ and $i_*$ are 
the obvious ones. It seems to be a folklore result that $i^*$ is
left adjoint to $i_*$ but we could not find a proof in
the literature so we provide one here; we also show
that $i_*$ is faithful, that its essential image is closed under
subquotients, and that it has
a right adjoint (this is what we mean by a ``closed immersion'').
Several people have been aware for some time that this is the
appropriate intuitive picture but, as far as I know,
no formal definition has been given and so no explicit proof has 
been given.

If $A$ is a graded subalgebra of $B$, 
the commutative results suggests there is 
should be a closed subspace $Z$ of $Y=\Projnc B$ and an affine map 
$g:Y \backslash Z \to \Projnc A$. Theorem
\ref{thm.ratl.map} establishes such a result under reasonable
hypotheses on $A$ and $B$. In fact, that result is set in a more
general context, namely a homomorphism $\phi:A \to B$ of graded
rings. Corollary \ref{cor.vdB} then says that if $\phi:A \to B$
and $B$ is a finitely presented left $A$-module, then there is an
affine map $g: \Projnc B \to \Projnc A$. This is a (special case of
a) non-commutative analogue of the commutative result that a finite
morphism is affine.

If $A$ is a quotient of a commutative polynomial ring, and
$A^{(n)}$ is the graded subring with components
$(A^{(n)})_i=A_{ni}$, then there is an isomorphism
of schemes $\Proj A \cong \Proj A^{(n)}$.  Verevkin \cite{V} 
proved that $\Projnc A \cong \Projnc A^{(n)}$ when $A$ is 
no longer commutative, but is connected and generated in degree one.
Theorem \ref{thm.Ver} shows that when $A$ is not required
to be generated in degree one, there is still a map  $\Projnc A \to
\Projnc A^{(n)}$, and Proposition \ref{prop.open.isom} describes open
subspaces on which this map is an isomorphism.

The results here are modelled on the commutative case, and none are
a surprise. In large part the point of this paper is to make the
appropriate definitions so that results from commutative algebraic
geometry carry over verbatim to the non-commutative setting.
Thus we formalize and make precise some of the terminology and
intuition in papers like \cite{BGK} and \cite{KKO}. 

In Example \ref{eg.wtd.Pn} we show how our results apply
to a quotient of a weighted polynomial ring to obtain a birational 
isomorphism $\Proj_{nc} A \to X=\Proj A$ where $X$ is a commutative subscheme of a
weighted projective space. It can happen that $X$ is singular whereas
$\Proj_{nc} A$ is smooth. Thus we can view
$\Projnc A \to \Proj A$ as a non-commutative resolution of singularities. 

We freely use basic notions and terminology for non-commutative
spaces from the papers \cite{Sm1}, \cite{Sm2}, and \cite{vdB}.

\medskip
{\bf Acknowledgements.}
This work was stimulated by Darin Stephenson's paper \cite{S}.
I would like to thank him for explaining his results, and
also suggesting that Proposition \ref{prop.split} should be true.
The final paragraph of this paper applies our results to Stephenson's
 algebras.
I thank James Zhang for bringing \cite{ADR} to my attention.

\section{Definitions and preliminaries}
\label{sect.defns}

Throughout this paper we assume that $A$ is a
{\sf locally finite $\NN$-graded algebra} over a field $k$. Thus
$A=A_0 \oplus A_1 \oplus \cdots$, and $\dim_k A_i < \infty$ for all
$i$. 
The {\sf augmentation ideal} $\fm$ of $A$ is $A_1 \oplus A_2 \oplus
\cdots$.
If $A_0$ is finite dimensional and $A$ is right
noetherian, then it follows that $\dim_k A_i < \infty$ for all 
$i$ because $A_{\ge i}/A_{\ge i+1}$ is a noetherian 
$A/\fm$-module. We write $\GrMod A$ for the category of
$\ZZ$-graded right $A$-modules and define
$$
\Tails A:=\GrMod A/\Fdim A,
$$
where $\Fdim A$ is the full subcategory 
consisting of direct limits of finite dimensional $A$-modules.
Equivalently, $\Fdim A$ consists of those modules every element of
which is annihilated by a suitably large power of $\fm$.
We write $\pi$ for the quotient functor $\GrMod A \to \Tails A$ and
$\omega$ for its right adjoint.

The {\sf projective space with homogeneous coordinate ring $A$} is
the space $X$ defined by $\Mod X:=\Tails A$. We write $\Projnc
A=(\Mod X,\cO_X)$ where $\cO_X$ denotes the image of $A$ in $\Tails
A$.

A {\sf closed subspace $Z$} of a space $X$ is a
full subcategory $\Mod Z$ of $\Mod X$ that is closed under
submodules and quotient modules in $\Mod X$ and such that the
inclusion functor $i_*:\Mod Z \to \Mod X$ has both a left adjoint
$i^*$ and a right adjoint $i^!$.
A map $\a:Y \to X$ of non-commutative spaces is a {\sf closed
immersion} if it is an isomorphism from $Y$ onto a closed subspace
of $X$.

The {\sf complement} $X\backslash Z$ to a closed subspace $Z$ is
defined by
$$
\Mod X \backslash Z := \Mod X /\sT,
$$
the quotient category of $\Mod X$ by the localizing subcategory
$\sT$ consisting of those $X$-modules $M$ that are the
direct limit of modules $N$ with the property that $N$ has a finite
filtration $N=N_n \supset N_{n-1} \supset \cdots N_1 \supset N_0=0$
such that each $N_i/N_{i-1}$ is in $\Mod Z$.
Because $\sT$ is a localizing category, there is an exact quotient
functor $j^*:\Mod X \to \Mod X \backslash Z$, and its right adjoint
$j_*:\Mod X \backslash Z \to \Mod X$. The pair $(j^*,j_*)$ defines
a map $j:X
\backslash Z \to X$. We call it an {\sf open immersion}.

We sometimes write $\Mod_ZX$ for the category $\sT$ and call it the
category of $X$-modules {\sf supported on} $Z$.

\medskip

Following Rosenberg \cite[pg. 278]{R},
a map $f:Y \to X$ is {\sf affine} if $f_*$ is faithful and 
has both a left adjoint $f^*$ and a right adjoint $f^!$.
If $f_*$ is faithful the counit $\id_Y \to f^!f_*$ is monic
and the unit $f^*f_* \to \id_Y$ is epic.

\bigskip

{\bf Watt's Theorem for graded modules.}
Let $A$ and $B$ be $\ZZ$-graded $k$-algebras.
We recall the analogue of Watt's Theorem proved by Del Rio 
\cite[Proposition 3]{ADR} that describes the $k$-linear 
functors $\GrMod A \to \GrMod B$ that have a 
right adjoint.

A bigraded $A$-$B$-bimodule is an $A$-$B$-bimodule 
$$
M=\bigoplus_{(p,q) \in \ZZ^2} {}_pM_q
$$
such that $A_i.{}_pM_q.B_j \subset {}_{i+p}M_{q+j}$ for all
$i,j,p,q \in \ZZ$.  
Write $\otimes$ for $\otimes_k$. If $L$ is a graded right 
$A$-module we define 
$$
L \barotimes_A M= \bigoplus_{q \in \ZZ} (L \barotimes_A M)_q
$$
where $(L \barotimes_A M)_q$ is the image of 
$\bigoplus_{p} (L_{-p} \otimes {}_pM_q)$ under the
canonical map $ L \otimes M \to L \otimes_A M$. 
This gives $L \barotimes_A M$ the structure of a graded right
$B$-module; it is a $B$-module direct summand of the usual tensor
product $L \otimes_A M$.
If $N$ is a graded right $B$-module we define 
$\uuHom_B(M,N)$ to consist of those $B$-module
homomorphisms that vanish on all except a finite number of
${}_pM_*$ and send each ${}_pM_q$ to $N_q$.
This is made into a graded right $A$-module by declaring that
$$
\uuHom_B(M,N)_p=\Hom_{\Gr B}({}_{-p}M_*,N).
$$
Then 
\begin{equation}
\label{eq.tens.hom.adj}
\Hom_{\Gr B}(L\barotimes_A M,N) \cong \Hom_{\Gr A}(L,\uuHom_B(M,N)).
\end{equation}

\begin{theorem}
[Del Rio] 
\cite{ADR}
Let $A$ and $B$ be graded $k$-algebras, and $F:\GrMod A \to \GrMod B$
a $k$-linear functor having a right adjoint. Then
$F \cong - \barotimes_A M$ where $M$ is the bigraded $A$-$B$-bimodule
$$
M= \bigoplus_{p \in \ZZ} F(A(p))
$$
with homogeneous components ${}_pM_q= F(A(p))_q$.

If $F$ also commutes with the twists by degree, then $F$ is
given by tensoring with a graded $A$-$B$-bimodule, 
say $V=\oplus_n V_n$. The corresponding $M$ in this case is $M=\oplus
V(p)$ with ${}_pM_q=V(p)_q$.

The left $A$-action on $M$ is given by declaring that $x \in A_i$ acts
on ${}_pM_*$ as $F(\l_x)$, where $\l_x:A(p) \to A(p+i)$ denotes left
multiplication by $x$.
\end{theorem}

\section{Maps induced by graded ring homomorphisms}
\label{sect.ratl.map}

Throughout this section we assume that $A$ and $B$ are 
locally finite $\NN$-graded algebras over a field $k$.

We consider the problem of when a homomorphism $\phi:A \to B$ 
of graded rings induces a map $g:\Projnc B \to \Projnc A$ and, if it
does, how the properties of $g$ are determined by the properties of 
$\phi$.

Associated to $\phi$ is an adjoint triple $(f^*,f_*,f^!)$
of functors
between the categories of graded modules. Explicitly,
$f^*=-\otimes_A B$, $f_*=-\otimes_B B_A$ is the restriction map,
and $f^!=\oplus_{p \in \ZZ} \Hom_{\Gr B}(B(-p),-)$. 
We wish to establish conditions on $\phi$ which imply
that these functors factor through the quotient categories in the 
following diagrams:
$$
\begin{CD}
\GrMod B @>{f_*}>> \GrMod A
\\
@V{\pi'}VV @VV{\pi}V
\\
\Tails B @. \Tails A
\end{CD}
\qquad \qquad
\begin{CD}
\GrMod B @<{f^*,\,f^!}<< \GrMod A
\\
@V{\pi'}VV @VV{\pi}V
\\
\Tails B @. \Tails A
\end{CD}
$$

\begin{lemma}
\label{lem.quotient.adjoint}
Let $\sA$ and $\sB$ be Grothendieck categories with localizing
subcategories $\sS \subset \sA$ and $\sT \subset \sB$.
Let $\pi:\sA \to \sA/\sS$ and
$\pi':\sB \to \sB/\sT$ be the quotient functors, and let $\omega$
and $\omega'$ be their right adjoints. 
Consider the following diagram of functors:
$$
\begin{CD}
\sA @>{F}>> \sB
\\
@V{\pi}VV @VV{\pi'}V
\\
\sA/\sS @. \sB/\sT.
\end{CD}
$$
\begin{enumerate}
\item{}
If $F(\sS) \subset \sT$,
then there is a unique functor $G:\sA/\sS \to \sB/\sT$ such that
$\pi'F=G\pi$. 
\item{}
If $H:\sB \to \sA$ is a right adjoint to $F$, then
$\pi H \omega'$ is a right adjoint to $G$.
\item{}
If $H$ is a right adjoint to $F$ and $G'$ is a right adjoint to $G$, 
then $H(\sT) \subset \sS$ if and only if $G'\pi' \cong  \pi H$.
\end{enumerate}
\end{lemma}
\begin{pf}
(1)
The existence and uniqueness of a functor $G$ such that
$\pi'F=G\pi$ is due to Gabriel \cite[Coroll. 2, p. 368]{Gab}.

(2)
To show that $G$ has a right adjoint it suffices to show that it is
right exact and commutes with direct sums.
If $\cM_\l$ is a collection of objects in $\sA/\sS$, then each is
of the form $\cM_\l =\pi M_\l$ for some object $M_\l$ in $\sA$.
Both $\pi'$ and $F$ commute with direct sums because they have
right adjoints, so $G\pi$ commutes with direct sums; $\pi$ also
commutes with direct sums.
Therefore
$$
G(\oplus \cM_\l)=G(\oplus \pi M_\l) \cong G\pi(\oplus M_\l) \cong \oplus
G\pi M_\l=\oplus G \cM_\l.
$$
Thus $G$ commutes with direct sums.

To see that $G$ is right exact consider an exact sequence
\begin{equation}
\label{eq.ses1}
0 \to \cL \to \cM \to \cN \to 0
\end{equation}
in $\sA/\sS$. By Gabriel \cite[Coroll. 1, p. 368]{Gab},
(\ref{eq.ses1}) is obtained by applying $\pi$ to an exact sequence
$0 \to L \to M \to N \to 0$ in $\sA$. Both $\pi'$ and $F$ are right
exact because they have right adjoints, so
$\pi'FL \to \pi'FM\to \pi'FN \to 0$ is exact. In other words
$G \cL \to G\cM \to G\cN \to 0$ is exact.

Hence $G$ has a right adjoint, say $G'$. It follows that $\omega
G'$ is a right adjoint to $G\pi$. But $G\pi=\pi' F$ has $H \omega'$
as a right adjoint so $\omega G' \cong H \omega'$. Since
$\pi\omega \cong \id_{\sA/\sS}$, $G' \cong \pi H \omega'$. Since a
right adjoint is only determined up to natural equivalence 
$\pi H \omega'$ is a right adjoint to $G$.

(3)
If $H(\sT) \subset \sS$, then $\pi H$ vanishes on $\sT$ so, by Gabriel
\cite[Corollaire 2, page 368]{Gab},
there is a functor $V:\sB/\sT \to \sA/\sS$ such that $V\pi'=\pi H$.
Thus $V \cong \pi H \omega'$, and this is isomorphic to $G'$ by
(2). Hence $G'\pi' \cong \pi H$.
Conversely, if $G'\pi' \cong \pi H$, then
$\pi H(\sT)=0$, so $H(\sT) \subset \sS$.
\end{pf}

{\bf Warning.}
The functor $H$ in part (2) of Lemma \ref{lem.quotient.adjoint} need not
have the property that $H(\sT)$ is contained in $\sS$. 
An explicit example of this is provided by taking
$\sB=\sA$, $F=H=\id_\sA$, $\sS=0$, and $\sT=\sB$.

\begin{theorem}
\label{thm.closed.proj}
Let $J$ be a graded ideal in an $\NN$-graded $k$-algebra $A$.
Then the homomorphism $A \to A/J$ induces a closed immersion 
$i:\Projnc A/J \to \Projnc A$. 
\end{theorem}
\begin{pf}
Write $X=\Projnc A$ and $Z=\Projnc A/J$.
Write $\fm=A_1 \oplus A_2 \oplus \cdots$. Thus $\Mod X=\GrMod
A/\Fdim A$.
We write $\pi:\GrMod A \to \Mod X$ for the quotient functor and
$\omega$ for a right adjoint to it. Similarly, $\pi':\GrMod A/J \to
\Mod Z$ is the quotient functor, and $\omega'$ is a right adjoint
to it. See \cite{V} and \cite[Section 2]{AZ} for more information
about this.

Let $f_*:\GrMod A/J \to \GrMod A$ be the inclusion functor.
It has a left adjoint $f^*=-\otimes_A A/J$, and a right adjoint
$f^!$ that sends a graded $A$-module to the largest submodule of it
that is annihilated by $J$.

By \cite[Coroll. 2, p. 368]{Gab} (=part(1) of Lemma
\ref{lem.quotient.adjoint}),
there is a unique functor $i_*$ such that 
$$
\begin{CD}
\GrMod A/J @>{f_*}>> \GrMod A
\\
@V{\pi'}VV @VV{\pi}V
\\
\Mod Z @>>{i_*}> \Mod X
\end{CD}
$$
commutes, and $i_*$ is exact because $f_*$ is \cite[Coroll. 3, p.
369]{Gab}. Thus
$i_*\pi'=\pi f_*$.

By Lemma \ref{lem.quotient.adjoint}(2),
$i_*$ has a right adjoint, namely $i^!:=\pi' f^! \omega$.
It is clear that $f^!$ sends $\Fdim A$ to $\Fdim A/J$, so 
$\pi' f^! \cong i^! \pi$ by  Lemma \ref{lem.quotient.adjoint}(3).

It is clear that $f^*$ sends $\Fdim A$ to $\Fdim A/J$, 
so by \cite[Coroll. 2, p. 368]{Gab}, there is a functor $i^*:\Mod X
\to \Mod Z$ such that $\pi' f^* = i^* \pi$. Since $f_*$ is right
adjoint to $f^*$, it follows from Lemma
\ref{lem.quotient.adjoint}(2) that
$\pi f_* \omega'$ is a right adjoint to $i^*$. But $\pi f_*\omega'=i_*
\pi'\omega' \cong i_*$. Hence $i^*$ is left adjoint to $i_*$.

We now show that $i_*$ is faithful. Since $i_*$ has a left and a
right adjoint it is exact, so it
suffices to show that if $i_*\cM=0$, then $\cM=0$.
Suppose that $i_*\pi'M=0$ for some
$M \in \GrMod A/J$.
Then $\pi f_* M=0$, and we conclude that $M$ is in $\Fdim A$, and hence
in $\Fdim A/J$; therefore $\pi' M=0$.
Hence $i_*$ is faithful.

We will show that $i_*$ is full after establishing the following
fact.

\underline{Claim:} $\omega\pi f_* \cong f_*\omega'\pi'$.
\underline{Proof:} Let $M \in \GrMod A/J$, let $\tau M$ denote
the largest submodule of $M$ that is in $\Fdim A/J$ (equivalently,
in $\Fdim A$), and set $\Mol=M/\tau M$. Then $\pi' M = \pi' \Mol$
and $\pi f_*M = \pi f_* \Mol$, so the two
functors take the same value on $M$ if and only if they take the
same value on $\Mol$. Hence we can, and will, assume that $M=\Mol$;
i.e., $\tau M=0$.

We must show that $\omega \pi M = \omega'\pi'M$.
By definition $\omega\pi M$ is the largest essential extension
$0 \to M \to \omega \pi M \to T \to 0$ such that $T \in \Fdim A$.
The definition of $\omega'\pi' M$ is analogous (although $T$ is
then required to belong to $\Fdim A/J$), so it suffices to
prove that $\omega \pi M$ is in $\GrMod A/J$. The surjective map
$\omega \pi M \otimes_A J \to (\omega \pi M)J$ is such that the
composition $M \otimes_A J \to \omega \pi M \otimes_A J \to (\omega
\pi M)J$ is zero, so there is a surjective map $T \otimes_A J \to
(\omega \pi M)J$. However, $T \otimes_A J$  belongs to $\Fdim A$
because $T$ does, so $(\omega \pi M)J \in \Fdim A$. This
implies that $M \cap (\omega \pi M)J \in \Fdim A$. But $\tau M=0$,
so $M \cap (\omega \pi M)J =0$, and it follows that $(\omega \pi
M)J=0$ because $M$ is essential in $\omega\pi M$. In other
words, $\omega\pi M \in \GrMod A/J$.
This completes the proof of the Claim.
$\square$

We have $f^*f_* \cong \id_{\GrMod A/J}$ and $\pi'\omega' \cong
\id_{\Mod Z}$, so
$$
i^*i_* \cong (\pi' f^* \omega)(\pi f_* \omega') \cong
\pi'f^*f_*\omega'\pi'\omega' \cong \id_{\Mod Z}.
$$
It follows from this that $i_*$ is full.

To see that $i_*$ is a closed immersion, 
it remains to check that $i_*(\Mod Z)$ is closed under submodules
and quotients in $\Mod X$.
Let $\cM \in \Mod Z$ and suppose that $0 \to \cL \to i_*\cM \to \cN
\to 0$ is an exact sequence in $\Mod X$. There is an exact sequence 
$0 \to \omega \cL \to \omega i_*\cM \to \omega \cN \to R^1\omega \cL$ 
in $\GrMod A$. 
Let $N$ denote the image of $\omega i_*\cM$ in $\omega \cN$.
Then $\cL \cong \pi\omega \cL$ and $\cN \cong \pi N$ because $\pi$
is exact and $R^1\omega\cL \in \Fdim A$.
Now $\cM=\pi' M$ for some $M \in \GrMod A/J$, so 
$\omega i_* \cM = \omega i_* \pi' M = \omega \pi f_* M \cong 
f_*\omega'\pi' M$ from which we conclude that $\omega i_* \cM$ is
annihilated by $J$. Therefore $\omega \cL$ is also annihilated by
$J$, so is of the form $f_*L$ for some $L \in \GrMod A/J$; hence
$\cL \cong \pi\omega \cL \cong \pi f_*L \cong i_* \pi'L \in
i_*(\Mod Z)$. Since $N$ is a quotient of $\omega i_*\cM$ it is also
annihilated by $J$, and a similar argument shows that $\cN \in
i_*(\Mod Z)$.

Since $f^*A=A/J$, $i^*\cO_X =\cO_Z$.
\end{pf}

We retain the notation of the theorem.

Because $i_*$ is fully faithful, we often view $\Mod Z$ as a full
subcategory of $\Mod X$ and speak of $Z=\Projnc A/J$ as a 
closed subspace of $X=\Projnc A$ and
call it the {\sf zero locus} of $J$.

It is {\it not} the case that every closed subspace of $\Projnc A$ 
is the zero locus of a two-sided ideal in $A$.
For example, if $A=k_q[x,y]$ is the ring defined by the relation
$yx=qxy$ where $0 \ne q \in k$, then $\Projnc A \cong \PP^1$, but the
closed points of $\PP^1$ are not cut out by two-sided ideals when
$q \ne 1$: for example, $(\a x+\b y)A$ is not a two-sided ideal when
$q \ne 1$ and $\a\b \ne 0$.
This is essentially due to the fact that the auto-equivalence $\cM \to
\cM(1)$ of $\Mod X$ induced by the degree shift on $A$ does not
generally send $Z$-modules to $Z$-modules. 

A more difficult question is whether every closed subspace of
$\Projnc A$ is the zero locus of a two-sided ideal in {\it some}
homogeneous coordinate ring of $\Projnc A$. We do not know the
answer to this question.

\begin{theorem}
\label{thm.ratl.map}
Suppose that $\phi:A \to B$ is a map of locally finite $\NN$-graded 
$k$-algebras. 
Write $X=\Projnc A$ and $Y=\Projnc B$.
Let $\fm$ be the augmentation ideal of $A$, and
let $I$ be the largest two-sided ideal of $B$ contained in
$\phi(\fm)B$. Let $Z\subset Y$ be the zero locus of $I$. 
If $B\phi(\fm)^n \subset \phi(\fm)B$ for some integer $n$, then
$\phi$ induces an affine map
$$
g:Y \backslash Z \to X.
$$
\end{theorem}
\begin{pf}
The category of modules over $Y \backslash Z$ 
is $\Mod Y/\Mod_Z Y$. This is equivalent to the quotient 
category $\GrMod B/\sT$ where $\sT$ consists of 
those modules $M$ with the property
that every element of $M$ is killed by some power of $I$.
Let $\pi':\GrMod B \to \GrMod B/\sT$ be
the quotient functor. We have functors $(f^*,f_*,f^!)$ between the
graded module categories and a diagram
$$
\begin{CD}
\GrMod B @>{f_*}>> \GrMod A
\\
@V{\pi'}VV @VV{\pi}V
\\
\Mod \, Y \backslash Z @. \Mod X.
\end{CD}
$$

To check that $f^*$ sends $\Fdim A$ to $\sT$ it suffices to check
that $f^*(A/\fm) \in \sT$ because $f^*$ commutes with direct limits
and with the degree twist $(1)$. However, $f^*(A/\fm)=B/\phi(\fm)B$
is in $\sT$ because $I \subset \phi(\fm)B$.
Hence
there is a unique functor $g^*:\Tails A \to (\GrMod B)/\sT$
satisfying $g^*\pi=\pi' f^*$.

To check that $f_*$ sends $\sT$ to $\Fdim A$
it suffices to check that $f_*(B/I)$ is in $\Fdim A$. 
However, $(B/I).\fm^n=B\phi(\fm)^n+I/I$; the hypothesis that
$B\phi(\fm)^n \subset \phi(\fm)B$ ensures that $B\phi(\fm)^n
\subset I$, so $(B/I).\fm^n=0$.
Hence there is an exact functor $g_*:\GrMod B/\sT \to
\Tails A$ such that $g_* \pi' = \pi f_*$. 

By Lemma \ref{lem.quotient.adjoint}(2), $g_*$ has a right adjoint
$g^!=\pi' f^! \omega$.

To show that $g_*$ is faithful we must show that if $M$ is a graded
$B$-module such that $g_*\pi' M=0$, then $M \in \sT$. Since $\pi f_*M
=0$, as an $A$-module $M$ is a direct limit $\liminj M_\l$ where 
each $M_\l$ is a finite dimensional $A$-module.
There is an epimorphism
$$
\liminj (M_\l \otimes_A B) \cong 
(\liminj M_\l) \otimes_A B \cong M \otimes_A B \to M
$$
of graded $B$-modules. 
Since $M_\l \otimes_A B$ equals $f^*M_\l$, it is in $\sT$; but $\sT$ is
closed under direct limits and quotients, so $M$ is in $\sT$.
Thus $g_*$ is faithful.
\end{pf}

The following consequence of the theorem slightly extends a result
of Van den Bergh  \cite[Proposition 3.9.11]{vdB}.

\begin{corollary}
\label{cor.vdB}
Let $\phi:A \to B$ be a homomorphism of graded rings such that
$B$ becomes a finitely presented graded left $A$-module. Then
$\phi$ induces an affine map $g:\Projnc B  \to  \Projnc A$.
\end{corollary}
\begin{pf}
If we apply $A/\fm \otimes_A -$ to a finite presentation of $B$
as a left $A$-module, we see that $B/\phi(\fm)B$ has finite
dimension. Thus, as a right $A$-module $B/\phi(\fm)B$ is annihilated 
by $\fm^n$ for some $n \gg 0$. Equivalently, $B \phi(\fm)^n  
\subset \phi(\fm)B$.
Thus the hypotheses of the theorem are satisfied. It remains to
show that $Z$ is empty.

Let $I$ denote the right annihilator in $B$ of $B/\phi(\fm)B$.
We have already observed that $\phi(\fm)^n \subset I$.
Since $A/\fm^n$ is finite dimensional, so is $A/\fm^n \otimes_A
B \cong B/\phi(\fm)^n B$. Thus $B/I$ is finite dimensional. Hence 
the zero locus of $I$ in $\Projnc B$ is empty.
\end{pf}

{\bf Remark.}
If, in Theorem \ref{thm.ratl.map}, $B_A$ is finitely presented,
then we have the  useful technical fact that  $g^!\pi=\pi'f^!$.
This follows from Lemma \ref{lem.quotient.adjoint}(3) once we show
that $f^!$ sends $\Fdim A$ to $\Fdim B$.
Let $M=\liminj M_\l$ be a direct limit of finite dimensional
$A$-modules. 
If $B$ is a finitely presented right $A$-module, then $\Hom_{\Gr
A}(B,-)$ commutes with direct limits, so $\Hom_{\Gr A}(N,\liminj
M_\l) = \liminj \Hom_{\Gr A}(B,M_\l)$; this is a direct limit of
finite dimensional $B$-modules because $B_A$ is a finitely
generated. Hence $f^!(\Fdim A) \subset \Fdim B$.

\section{The Veronese mapping}
\label{sect.veronese}

Throughout this section $A$ is a locally finite $\NN$-graded 
$k$-algebra and $n$ is a positive integer. 
The {\sf $n^{\th}$ Veronese} subalgebra $A^{(n)}$ is defined by
$$
A^{(n)}_i:=A_{ni}.
$$
It is a classical result in algebraic geometry that if $A$ is 
a finitely generated commutative connected graded
$k$-algebra generated in degree one, then $\Proj A \cong \Proj
A^{(n)}$. This isomorphism is implemented by the Veronese
embedding.

Verevkin proved a non-commutative version of this result when $A$ is
noetherian and generated in degree one \cite[Theorem 4.4]{V}. 

Theorem \ref{thm.Ver} and Proposition \ref{prop.open.isom}
show what happens when $A$ need not be commutative
and need not be generated in degree one.

\begin{theorem}
\label{thm.Ver}
Let $A$ be a right noetherian locally finite $\NN$-graded $k$-algebra.
Fix a positive integer $n$.
There is a map $g:\Projnc A\to \Projnc A^{(n)}$.
Furthermore, $g_*$ has a right adjoint.
\end{theorem}

We will use the notation $X:=\Projnc A^{(n)}$ and $X':=\Projnc A$.

We need two preliminary results before proving the theorem. 
First we explain how
the functors defining the map $g:\Projnc A  \to \Projnc A^{(n)}$ in the theorem are
induced by functors between the categories $\GrMod A$ and 
$\GrMod A^{(n)}$.

If $L$ is a graded $A$-module we define the graded $A^{(n)}$-module 
$L^{(n)}$ by
$$
L^{(n)}_i:=L_{ni}.
$$
The rule $L \mapsto L^{(n)}$ extends to give an exact functor
$$
f_*: \GrMod A \to \GrMod A^{(n)}.
$$
The functor $f_*$ is not faithful when $n \ge 2$ because
$f_*((A/\fm)(1))=0$.

The next result shows that $f_*$ has a left
adjoint $f^*$ and a right adjoint $f^!$. 

\begin{proposition}
\label{prop.Ver}
Let $A$ be a locally finite $\NN$-graded $k$-algebra.
Fix a positive integer $n$.
Let $W$ and $W'$ be the spaces with module categories
$$
\Mod W=\GrMod A^{(n)}
$$
and
$$
\Mod W'=\GrMod A.
$$
Then there is a map $f:W' \to W$ with direct image functor
given by $f_*L=L^{(n)}$.
\end{proposition}
\begin{pf}
It is clear that  $f_*$ is an exact functor commuting with direct sums.
By the graded version of Watt's Theorem, $f_* \cong -\barotimes_A M$
where
$$
M= \bigoplus_{p \in \ZZ} A(p)^{(n)}
$$
with components ${}_pM_q=(A(p)^{(n)})_q=A(p)_{nq}.$
The right action of $A^{(n)}$ is given by right multiplication, and
each $A(p)^{(n)}$ is a right $A^{(n)}$-submodule. The left action
of $A$ is given by left multiplication whereby $a \in A_i$ acts by
sending $A(p)_{nq}$ to $A(p+i)_{nq}$.

Define $f^*:\GrMod A^{(n)} \to \GrMod A$ by
$f^*N=N \otimes_{A^{(n)}} A$ with the usual right
action of $A$, and grading given by 
$$
(N \otimes_{A^{(n)}} A)_s= \sum_{ni+j=s} N_i \otimes A_j.
$$
It is not hard to show that $f^*$ is a left adjoint to $f_*$. 
Therefore $f^* \cong - \barotimes_{A^{(n)}} Q$ where
$$
Q=\bigoplus_{p \in \ZZ} f^*(A^{(n)}(p)) \cong \bigoplus_{p \in \ZZ}
A(np);
$$
multiplication $A^{(n)}(p) \otimes_{A^{(n)}} A \to A(np)$
gives an isomorphism of graded right $A$-modules.
Thus ${}_pQ_*\cong A(np)$ with its usual grading.
One can verify directly
that $f_* \cong \uuHom_A(Q,-)$.

The right adjoint to $f_*$ is the functor $f^!=\uuHom_{A^{(n)}}(M,-)$.
If $N$ is a graded right $A^{(n)}$-module
$$
(f^!N)_i=\Hom_{\Gr A^{(n)}} ({}_{-i}M_*,N) = \Hom_{\Gr A^{(n)}}
(A(-i),N).
$$

If $N$ is a graded $A^{(n)}$-module, then $f_*f^*(N)=N$ so $f_*f^*$ is
naturally equivalent to $\id_X$. 
\end{pf}

Let $\pi': \GrMod A \to \Tails A$ and $\pi:\GrMod A^{(n)}
\to \Tails   A^{(n)}$ be the quotient functors. 
To prove Theorem \ref{thm.Ver}, we must
find functors $g^*$, $g_*$,
and $g^!$ making the following diagrams commute:
$$
\begin{CD}
\GrMod A @>{f_*}>> \GrMod A^{(n)}
\\
@V{\pi'}VV @VV{\pi}V
\\
\Tails A @>>{g_*}> \Tails  A^{(n)}.
\end{CD}
\qquad
\qquad
\begin{CD}
\GrMod A @<{f^*,f^!}<< \GrMod A^{(n)}
\\
@V{\pi'}VV @VV{\pi}V
\\
\Tails A @<<{g^*,g^!}< \Tails  A^{(n)}.
\end{CD}
$$

Since $f_*$ sends $\Fdim A$ to $\Fdim A^{(n)}$ there
is a functor $g_*:\Tails A \to \Tails A^{(n)}$ such that $g_*\pi'=\pi
f_*$.

To ensure that $f^*$ and $f^!$ induce functors between the quotient 
categories we must impose a noetherian hypothesis. 
Although there is no noetherian hypothesis in Proposition 
\ref{prop.Ver}, in Theorem \ref{thm.Ver} it is assumed that $A$ is
right noetherian. This hypothesis ensures that $f^*$ sends $\Fdim
A^{(n)}$ to $\Fdim A$. 

\begin{lemma}
\label{lem.ver.noeth}
Let $A$ be a locally finite $\NN$-graded $k$-algebra.
\begin{enumerate}
\item{}
$f_*$ sends noetherian $A$-modules to noetherian $A^{(n)}$-modules;
\item{}
if $A$ is right noetherian so is $A^{(n)}$, and $A$ is a finitely
generated right $A^{(n)}$-module.
\item{}
if $A$ is left noetherian, then $f^*$ sends $\Fdim A^{(n)}$ to
$\Fdim A$;
\item{}
if $A$ is left and right noetherian, then there is a functor
$g^*:\Tails A^{(n)} \to \Tails A$ such that $g^*\pi=\pi'f^*$;
\item{}
if $A$ is left noetherian, then $f^!$ sends $\Fdim A^{(n)}$ to
$\Fdim A$;
\item{}
if $A$ is left and right noetherian, then there is a functor
$g^!:\Tails A^{(n)} \to \Tails A$ such that $g^!\pi=\pi'f^!$;
\end{enumerate}
\end{lemma}
\begin{pf}
(1)
Let $M$ be a right noetherian graded $A$-module. 
If $N$ is a submodule of $M^{(n)}$ then $N=NA \cap M^{(n)}$. Hence
any proper ascending chain of submodules in $M^{(n)}$ would give a
proper ascending chain of submodules of $M$ by multiplying
by $A$. Since $M$ contains no such chain, neither does $M^{(n)}.$

(2)
Applying (1) to $M=A$ shows that
$A^{(n)}$ is right and left noetherian.
If $N=A \oplus A(1) \oplus \cdots \oplus A(n-1)$, then $N^{(n)}
\cong A$ as a left and as a right $A^{(n)}$-module, so applying 
(1) to $N$ gives the rest of the result.

(3)
If $N$ is a finite dimensional right $A^{(n)}$-module, then $N
\otimes_{A^{(n)}} A$ is a finite dimensional $A$-module because $A$
is a finitely generated {\it left} $A^{(n)}$-module.

(4)
This follows at once from (3).

(5)
We must show that if $N$ is a finite dimensional graded right
$A^{(n)}$-module, then $\uuHom_{A^{(n)}}(M,N)$ is finite
dimensional, where
$$
M= \bigoplus_{p \in \ZZ} A(p)^{(n)}.
$$
The degree $-p$ component of $\uuHom_{A^{(n)}}(M,N)$ is
$\Hom_{\Gr A^{(n)}}(A(p)^{(n)},N)$, and it suffices to show that this
is zero for almost all $p$ and is always finite dimensional.
By the noetherian hypothesis, $A(p)^{(n)}$ is a finitely generated
$A^{(n)}$-module, so $\Hom_{\Gr A^{(n)}}(A(p)^{(n)},N)$ has finite
dimension. We must show that $\Hom_{\Gr
A^{(n)}}(A(p)^{(n)},N)$ is zero if $|p|$ is sufficiently large.

Fix $p$. Now
$$
A(p+nj)^{(n)} \cong A(p)^{(n)}(j)
$$
so 
$$
\Hom_{\Gr A^{(n)}}(A(p+nj)^{(n)},N) \cong 
\Hom_{\Gr A^{(n)}}(A(p)^{(n)},N(-j)).
$$
Since $A(p)^{(n)}$ is finitely generated and $N$ is finite
dimensional, when $|j|$ is sufficiently large
 $\Hom_{\Gr A^{(n)}}(A(p)^{(n)},N(-j))$ is zero.

It therefore follows that 
$\Hom_{\Gr A^{(n)}}(A(p)^{(n)},N)$ is zero for $|p|$ sufficiently
large. 
This completes the proof of (5), and (6) follows from this.
\end{pf}

{\bf Proof of Theorem \ref{thm.Ver}.}
By Lemma \ref{lem.ver.noeth} there are functors $g^*,g_*,g^!$ 
between the categories $\Mod \Projnc A^{(n)}=\Tails A^{(n)}$ and $\Mod \Projnc A=\Tails A$ satisfying
$$
g^*\pi = \pi' f^*, \quad 
g_* \pi'= \pi f_*, \quad
g^! \pi = \pi' f^!.
$$
Applying Lemma \ref{lem.quotient.adjoint} to $f^*$ we see that
$g^*$ has $\pi f_* \omega'$ as a right adjoint. But this is
naturally equivalent to $g_*$, so $g_*$ is a right adjoint to $g^*$. 
Similarly, $g^!$ is a right adjoint to $g_*$.
Since $f^*A^{(n)}=A$, $g_*\cO_{\Projnc A^{(n)}}=\cO_{\Projnc A}$. 
This completes the proof of  Theorem \ref{thm.Ver}.
\hfill $\square$

In the next result $\Proj A$ is the usual commutative scheme viewed
as a non-commutative space with module category $\Qcoh(\Proj A)$.

\begin{corollary}
\label{cor.ver}
If $A$ is a finitely generated commutative connected graded
$k$-algebra, there is a map $g:\Projnc A \to \Proj A$.
Furthermore, $g_*$ has a right adjoint $g^!$.
\end{corollary}
\begin{pf}
For some sufficiently large $n$, 
$A^{(n)}$ is generated in degree one so 
$$
\Tails A^{(n)} \cong \Tails A \cong \Qcoh \Proj A.
$$
Therefore Theorem \ref{thm.Ver} gives the result.
\end{pf}

{\bf Remarks. 1.}
Since $g_*$ has both a left and a right adjoint it is exact, and
hence its right adjoint $g^!$ preserves injectives. There is
therefore a convergent spectral sequence
$$
\Ext^p_{\Projnc A}(M,R^q g^!N) \Rightarrow \Ext^{p+q}_{\Projnc
A^{(n)}}(g_*M,N)
$$
for $M$ and $N$ modules over $\Projnc A$ and $\Projnc A^{(n)}$
respectively.

{\bf 2.}
If $J$ is a two-sided ideal of $A$, then the natural map $A \to
A/J$ induces an isomorphism $A^{(n)}/J^{(n)} \to (A/J)^{(n)}$, so
there is a commutative diagram
$$
\begin{CD}
\Projnc A/J @>{i}>> \Projnc A
\\
@VVV @VV{g}V
\\
\Projnc A^{(n)}/J^{(n)} @>>> \Projnc A^{(n)}
\end{CD}
$$
where the horizontal maps are the natural closed immersions.

\begin{proposition}
\label{prop.birat.isom}
Let $A$ be a right noetherian, locally finite, $\NN$-graded
$k$-algebra. Suppose that $A$ is prime and $\Fract_{gr} A$ 
contains a copy of $A(n)$ for all $n \in \ZZ$. Then 
\begin{enumerate}
\item{}
$\Projnc A$ and $\Projnc A^{(n)}$ are integral spaces
in the sense of \cite{Sm2}, and
\item{}
$g:\Projnc A \to \Projnc A^{(n)}$ is a birational isomorphism
in the sense that it induces an isomorphism between the function
fields.
\end{enumerate}
\end{proposition}
\begin{pf}
That  $X=\Projnc A$ is an integral space is proved in \cite[Th.
4.5]{Sm2}. It is also shown there that the function field of $X$
is 
$$
\Fract_{gr} A:=\{ab^{-1} \; | \; a,b\in A \hbox{ homogeneous
of the same degree and $b$ is regular}\}.
$$
It is clear that $\Fract_{gr} A^{(n)} \subset \Fract_{gr} A$, and
the reverse inclusion follows from the observation that
$ab^{-1}=ab^{n-1}b^{-n}$.
\end{pf}

{\bf Remarks. 1.}
If $A$ is prime noetherian and has a regular element of degree $d$
for all $d \gg 0$, then  $\Fract_{gr} A$ contains a copy of $A(n)$
for all $n \in \ZZ$, so the previous result applies.

{\bf 2.}
If $z$ is a normal regular element, then the complement to the zero
locus of $z$ in $\Projnc A$ is $\Projnc A[z^{-1}]$; the category
of modules over this is $\GrMod A[z^{-1}]$; if $d$ is the
smallest positive integer such that $A[z^{-1}]$ has a unit of
degree $d$, this is an affine space with coordinate ring
$$
\begin{pmatrix}
R_0 & R_1 & \cdots & R_{d-1}
\\
R_{-1} & R_0 & \cdots & R_{d-2}
\\
\vdots &&& \vdots
\\
R_{-d+1} & R_{-d+2} & \cdots & R_0
\end{pmatrix}
$$
where $R=A[z^{-1}]$.

{\bf 3.}
In the previous result, if $s$ and $t$ are homogeneous
regular elements of relatively prime degrees in $A$, and $st$ is
normal, then $A[(st)^{-1}]$ has a unit of degree one, so the open
complement to the zero locus of $st$ in $\Projnc A$ is the affine
space with coordinate ring $A[(st)^{-1}]_0$. This ring is equal to
$A^{(n)}[(st)^{-n}]_0$, so the open complement is 
is isomorphic to the open complement to the zero locus of $(st)^n$
in $\Projnc A^{(n)}$.

{\bf 4.}
If $\Fract_{gr} A$ fails to contain a copy of $A(n)$ for all $n$, 
the map
$\Projnc A \to \Projnc A^{(n)}$ need not be a birational isomorphism. 
For example, take $A=k[x]$ with $\deg x=2$.

\begin{example}
\label{eg.bad}
If $A$ is not generated in degree one, then $g_*$ need not be faithful.

Let $A=k[x,z]$ be the polynomial ring with $\deg x=1$ and 
$\deg z=n\ge 2$. The image under $\pi$ of $M=A/(x)$ is a simple module 
$\cO_p$ in $\Projnc A$. We have $\cO_p(1) \ne 0$, but $(M(1))^{(n)}=0$, 
so $g_*(\cO_p(1))=0$. 
\end{example}

One might anticipate that $g:\Projnc A \to \Projnc A^{(n)}$ is 
an isomorphism on suitable open subspaces: in the previous 
example, $g$
restricts to an isomorphism from the complement to
the zero locus of $x$ in $\Projnc A$ to the complement to the zero 
locus of $x$ in $\Projnc A^{(n)}$.
We prove a general result of this type in Proposition
\ref{prop.open.isom}. First we need a lemma.

\smallskip

For each integer $r$, define 
$$
A^{(n)+r}: =  \sum_{j \in \ZZ} A_{nj+r}.
$$
Obviously, $A^{(n)+r}A^{(n)+s} \subset A^{(n)+r+s}$ so each
$A^{(n)+r}$ is an  $A^{(n)}$-$A^{(n)}$-bimodule, and these bimodules 
depend only on $r$(mod$\,n$).
Define
$$
I_r := A^{(n)+r}A=\sum_{j \in \ZZ} A_{nj+r}A
$$
and
$$
I:=\bigcap_{r \in \ZZ} I_r= I_1 \cap I_2 \cap \cdots \cap I_n.
$$
Although $I_r$ is in general only a right ideal of $A$, $I_r^{(n)}$
is a two-sided ideal of $A^{(n)}$.
Since $A_qI_r \subset I_{q+r}$, $I$
is a two-sided ideal of $A$. 

Notice that $A^{(n)+r}A^{(n)-r} = I_r^{(n)}$. 

\begin{lemma}
\label{lem.I(n)}
With the above notation, $I^{2n} \subset I^{(n)}A.$
\end{lemma}
\begin{pf}
From the containment
$$
I^2 \; \subset \; I_rI = A^{(n)+r}I \;  \subset \;  A^{(n)+r} I_{n-r}
 =  A^{(n)+r} A^{(n)+n-r}A =I_r^{(n)}A,
$$
it follows that 
$$
I^{2n} \; \subset \; I_1^{(n)} I^{2n-2}\; \subset \; 
I_1^{(n)}I_2^{(n)} I^{2n-4}
\;\subset \cdots \;\subset \; I_1^{(n)} \cdots I_n^{(n)}A.
$$
But this last term is contained in 
$$
\bigl(I_1^{(n)} \cap \cdots \cap I_n^{(n)}\bigr)A =I^{(n)} A,
$$
which completes the proof.
\end{pf}

\begin{proposition}
\label{prop.open.isom}
With the hypotheses of Theorem \ref{thm.Ver}, the map $g$ restricts
to an isomorphism $g:\Projnc A \backslash Z'  \to 
\Projnc A^{(n)} \backslash Z$
where $Z'$ and $Z$ are the zero loci of $I$ and $I^{(n)}$
respectively. 
\end{proposition}
\begin{pf}
 Write $X'=\Projnc A$, $X=\Projnc A^{(n)}$,
$U'=X'\backslash Z'$ and $ U=X \backslash Z$.
Write $\a:U \to X$ and $\b:U' \to X'$ for the inclusions.
We will use Lemma \ref{lem.quotient.adjoint} to show there is
an isomorphism $h:U' \to U$ such that the diagram
$$
\begin{CD}
U'@>{\b}>> X' = \Projnc A
\\
@V{h}VV @VV{g}V
\\
U @>>{\a}> X =\Projnc A^{(n)} 
\end{CD}
$$
commutes.

Let $\sT$ be the localizing subcategory of $\GrMod A$ consisting of
those modules $L$ such that every element of $L$ is killed by a
suitably large power of $I$. Let $\sS$  be the localizing
subcategory of $\GrMod A^{(n)}$ consisting of those modules $N$
such that every element of $N$ is killed by a suitably large power
of $I^{(n)}$. These two localizing subcategories contain all the
finite dimensional modules. The spaces $U$ and $U'$
are defined by
$$
\Mod U:= (\GrMod A^{(n)})/\sS 
\qquad \hbox{and} \qquad 
\Mod U':= (\GrMod A)/\sT.
$$

Let $f_*$ and $f^*$ be the functors defined in Proposition
\ref{prop.Ver}. We will
show that $f_*(\sT) \subset \sS$ and $f^*(\sS) \subset \sT$.
The first of these inclusions is obvious: if every
element of an $A$-module $L$ is annihilated by a power of $I$, then
every element of $L^{(n)}$ is annihilated by a power of $I^{(n)}$.
To show that $f^*(\sS) \subset \sT$ it suffices to show that
$f^*(A^{(n)}/I^{(n)})$ belongs to $\sT$. But
 $f^*(A^{(n)}/I^{(n)}) \cong A/I^{(n)}A$, and by Lemma
\ref{lem.I(n)}, $A/I^{(n)}A$ is annihilated by $I^{2n}$ so belongs
to $\sT$.

We now use Lemma \ref{lem.quotient.adjoint} in the context of the
following diagram
$$
\begin{CD}
\GrMod A^{(n)} @>{f^*}>> \GrMod A
\\
 @V{\a^*\pi}VV @VV{\b^{*}\pi'}V
\\
\Mod U @. \Mod U'.
\end{CD}
$$
Because $f^*(\sS) \subset \sT$, there exists a functor $h^*:\Mod U
\to \Mod U'$ such that $h^*\a^*\pi = \b^{*}\pi'f^*$.
Because $f_*$ is right adjoint to $f^*$,
$h_*:=\a^*\pi f_*\omega'\b_*$ is right adjoint to $h^*$. 
Thus, $h^*$ and $h_*$ define a map $h:U' \to U$.
Since $g_*\pi'=\pi f_*$, a computation gives 
$\a_*h_* \cong g_*\b_*$. Therefore $\a h = g\b$.

It remains to show that $h$ is an isomorphism.

The unit $\id_{\GrMod A^{(n)}} \to f_*f^*$ is an
isomorphism because the natural map $L \to (L \otimes_{A^{(n)}}
A)^{(n)}$ is an isomorphism for all $L \in \GrMod A^{(n)}$.
Because $f_*(\sT) \subset \sS$, part (3) of Lemma
\ref{lem.quotient.adjoint} gives $h_*\b^{*}\pi' \cong \a^*\pi f_*$.
Therefore,
$$
h_*h^*\cong h_*h^*\a^*\pi\omega \a_* =h_*\b^{*}\pi'f^* \omega \a_* 
\cong \a^*\pi f_* f^* \omega \a_* \cong \id_U.
$$

To show that the natural transformation $h^*h_* \to \id_{U'}$ is an
isomorphism, we first consider the natural transformation 
$f^*f_* \to \id_{\GrMod A}$. 
For an $A$-module $M$ this is the multiplication map
$$
f^*f_*M=M^{(n)} \otimes_{A^{(n)}} A \to M.
$$
We claim that the kernel and cokernel of this map belong to $\sT$.

Suppose that $\sum m_i \otimes a_i \in M^{(n)} \otimes_{A^{(n)}} A$
is in the kernel. Then $\sum_i m_i a_i =0$. 
By taking homogeneous components we can reduce to the case where 
each $a_i$ belongs to $A_{n-r}+A_{2n-r}+\cdots$ for some
$r \in \{1,\ldots,n\}$. Then, if $b \in A_{nj+r}$ for some $j$, then
$$
\biggl( \sum_i m_i \otimes a_i \biggr) b= \sum_i m_i \otimes a_ib =
\sum_i m_i  a_ib \otimes 1 =0.
$$
Thus $\bigl( \sum_i m_i \otimes a_i \bigr) I_r=0$. Hence the kernel
is annihilated by $I$, so belongs to $\sT$.
The cokernel of $f^*f_*M \to M$ is $M/M^{(n)}A$. 
If $r \in \{1,\ldots,n\}$, then
$M_{nj-r}I_r \subset M^{(n)}A$. Hence $I$ annihilates the cokernel.

Since the kernel and cokernel of $f^*f_* \to \id$ belong to $\sT$,
there is an isomorphism $\b^{*}\pi' f^*f_* \to \b^{*}\pi'$. Hence
$$
h^*h_* = h^* \a^*\pi f_*\omega'\b_* \cong \b^{*}\pi'f^* f_*\omega'\b_*
\cong \id_{U'}.
$$
\end{pf}

In Example \ref{eg.bad}, $I_1=(x)$ and $I_2=A$, so
$I=(x)$, whence $Z'$ is the zero locus of $x$.
This explains why we need to
remove the zero locus of $x$ to get the isomorphism.

If $A$ is generated in degree one, then $I_r=A_{\ge r}$ for $r \in
\{0,1,\ldots,n-1\}$, so $A/I_r \in \Fdim A$, whence $A/I \in \Fdim A$.
It follows that $Z$ and $Z'$ are empty, and therefore $U=X$ and
$U'=X'$. We therefore recover Verevkin's result $X \cong X'$ when
$A$ is generated in degree one over $A_0$.

\begin{example}
\label{eg.wtd.Pn}
Let $A$ be a weighted polynomial ring. That is,
$A=k[x_0,\ldots,x_n]$ where $\deg x_i=q_i \ge 1$. Write
$Q=(q_0,\ldots,q_n)$. Then $\PP^n_Q:=
\Proj A$ is called a {\sf weighted projective space}. It is
isomorphic to the quotient variety $\PP^n/\mu_Q$, where $\mu_Q =
\mu_{q_0} \times \cdots \times \mu_{q_n}$. There is a large integer
$d$ such that $A^{(d)}$ is generated in degree one. Hence
$$
\PP^n_Q = \Proj A = \Proj A^{(d)} \cong \Projnc A^{(d)}.
$$
By Theorem \ref{thm.Ver}, there is a map 
$$
g: \Projnc A \to \Projnc A^{(d)} \cong \PP^n_Q.
$$
This is an isomorphism on an open subspace by Proposition
\ref{prop.open.isom}. Since $A$ has
global homological dimension $n+1$, $\Projnc A$ has global
homological dimension $n$. We therefore think of  $\Projnc A$ as a
smooth space of dimension $n$ and the map $g$ as a non-commutative
resolution of $\PP^n_Q$. 
Let $X \subset \PP^n_Q$ be the closed subscheme cut out by an 
ideal $J$ in $A$.  Then there is a commutative diagram 
$$
\begin{CD}
\Projnc A/J @>{i}>> \Projnc A
\\
@V{f}VV @VV{g}V
\\
X @>>> \PP^n_Q
\end{CD}
$$
in which $f$ is a birational isomorphism and $i$ is a closed immersion.
It can happen that $\Projnc A/J$ is smooth even when $X$ is singular.
Thus $\Projnc A$ is a ``non-commutative resolution'' of 
$X$. An interesting case to examine in some detail is that
where $X$ is an orbifold of a Calabi-Yau three-fold.
\end{example}

If $A=k[x]$ with $\deg x=2$, and $n=2$, then $\Projnc A \cong \Spec k
\times k$ and $\Proj A \cong \Spec k$. Furthermore, $Z'=X'$ and $Z=X$.
This is a special case of the next result, the truth of which was
suggested by Darin Stephenson.

\begin{proposition}
\label{prop.split}
Let $A$ be a locally finite $\NN$-graded $k$-algebra such that
$A_i=0$ whenever $i \not\equiv 0$(mod $n$).
Then $\Projnc A$ is isomorphic to the disjoint union of 
$n$ copies of $\Projnc A^{(n)}$.
\end{proposition}
\begin{pf}
Let $p_r:\GrMod A \to \GrMod A$ be the functor defined by
$$
p_r(M)=\bigoplus_{i \in \ZZ} M_{r+in}
$$
on objects, and $p_r(\theta) = \theta|_{p_r(M)}$ whenever $\theta
\in \Hom_{\Gr A}(M,N)$.
The hypothesis on $A$ ensures that each $p_r(M)$ is a graded
$A$-submodule of $M$, so $p_r$ is indeed a functor from $\GrMod A$
to itself. 
It is clear that $\id_{\GrMod A} = p_0 \oplus \cdots \oplus
p_{n-1}$, where this direct sum is taken in the abelian category of
$k$-linear functors from $\GrMod A$ to itself; essentially, this is
the observation that $M=p_0(M) \oplus \cdots \oplus
p_{n-1}(M)$, and that any map $\theta:M \to N$ of graded
$A$-modules respects this decomposition.
Furthermore, each $p_r$ is idempotent and the $p_r$s are mutually
orthogonal. It follows from this that there is a decomposition of
$\GrMod A$ as a product of categories, each component being the
full subcategory on which $p_r$ is the identity.

It is clear that the shift functor $(1)$ cyclicly permutes these
subcategories, so they are all equivalent to one another and $(n)$
is an auto-equivalence of each component.
However, any one of these categories together with its
autoequivalence $(n)$ is equivalent to $\GrMod A^{(n)}$ with its
auto-equivalence $(1)$. 
Thus $\GrMod A$ is equivalent to the product of $n$ copies of
$\GrMod A^{(n)}$.

Finally, this decomposition descends to the Tails categories.
\end{pf}

\section{An Ore extension and an example}

The morphism
\begin{align*}
p:\PP^n \backslash \{(0,\ldots,0,1)\} & \to \PP^{n-1}
\\
(\a_0,\ldots,\a_n) & \mapsto  (\a_0,\ldots,\a_{n-1})
\end{align*}
is called the {\sf projection with center} $(0,\ldots,0,1)$. 
This section examines a non-commutative analogue of this
basic operation.

Consider a connected graded $k$-algebra $R$ and a connected 
graded Ore extension
$$
S=R[t;\s,\d]
$$
with respect to a graded automorphism $\s$ and a graded
$\s$-derivation $\d$ of degree $n \ge 1$. Thus $S=\oplus_{n=0}^\infty
Rt^n$ and $tr=r^\s t+\d(r)$ for all $r \in R$. Since $\d(R_i) \subset
R_{i+n}$ for all $i$, by setting $\deg t=n$, $S$ becomes a 
connected graded algebra. 

One expects that the inclusion map $R \to S$ induces a map $\Projnc
S \to \Projnc R$. Indeed, the
projection map above can be obtained as a special case of this.

Let $\fm$ denote the augmentation ideal of $R$.
Since $\d(\fm) \subset \fm$, $\fm S$ is a two-sided ideal of $S$.
Furthermore, $S/\fm S \cong k[t]$ as graded rings.

\begin{proposition}
With the above notation, let $Z$ denote the zero locus of 
$\fm S$ in $Y$.
\begin{enumerate}
\item{}
$Z \cong \Spec k^{\times n}$.
\item{}
There is an affine map $g:\Projnc S \backslash Z \to \Projnc R$.
\end{enumerate}
\end{proposition}
\begin{pf}
(2)
The existence of $g$ is a special case of Theorem
\ref{thm.ratl.map}. That proposition applies because $S\fm \subset
\fm S$.

(1)
The quotient ring $S/\fm S$ is isomorphic to the polynomial ring
$k[t]$ with $\deg t=n$, so this follows from Proposition
\ref{prop.split}. 
\end{pf}

We think of $\Projnc S$ as a ``cone over $\Projnc R$ with vertex $Z$''.
It would be interesting to describe the ``fibers'' of the map $g$.

\medskip

When $\deg t >1$, the Ore extension $S=R[t;\s,\d]$ is not generated 
by its elements of degree one. This sometimes causes technical
problems; however, if $R$ is generated in degree one, then the 
$n^{\th}$ Veronese $S^{(n)}$ is generated in degree one.
We can then combine Theorems \ref{thm.ratl.map} and \ref{thm.Ver} 
to analyze the space with homogeneous coordinate ring $S$ as
follows.

\begin{proposition}
\label{prop.cones}
The inclusion of the $n$-Veronese subalgebras of $S$ and $S/\fm S$
gives a commutative diagram of rings and an induced commutative
diagram of spaces as in the following diagram:
$$
\begin{CD}
k[t]^{(n)} @>>> k[t]=S/S\fm
\\
@AAA @AAA
\\
S^{(n)} @>>> S
\\
@AAA @AAA
\\
R^{(n)} @>>> R
\end{CD}
\hskip 1in
\begin{CD}
\Spec k \cong v @<<< Z' \cong \Spec k^{\times n}
\\
@VVV @VVV
\\
\Projnc S^{(n)}  @<{g}<< \Projnc S 
\\
@AAA @AAA
\\
\Projnc S^{(n)} \backslash\{v\} @<<< \Projnc S \backslash Z'
\\
@V{\a}VV @VV{\b}V
\\
\Projnc R^{(n)} @<<< \Projnc R
\\
\end{CD}
$$
\end{proposition}

{\bf An application.}
In \cite{S}, a family of three-dimensional Artin-Schelter 
regular algebras $A$ is constructed and studied. Although the
algebraic properties of $A$ are quite well understood, our
understanding of the corresponding geometric object $\Projnc A$ is
rudimentary. The algebras are of the form
$A=R[t;\s,\d]$ with $\deg t=2$ and $R$ a two-dimensional
Artin-Schelter regular algebra generated in degree one.
It is well-known that $R$ and its Veronese subalgebras are (not
necessarily commutative) homogeneous coordinate rings of $\PP^1$.
By Proposition \ref{prop.cones}, there is a commutative diagram of
spaces and maps
\begin{equation}
\label{eq.maps}
\begin{CD}
\Spec k @<<< Z' \cong \Spec k^{\times n}
\\
@VVV @VVV
\\
\Projnc A^{(n)} @<{g}<< \Projnc A
\\
@AAA @AAA
\\
\Projnc A^{(n)} \backslash \{v\} @<{g}<< \Projnc A \backslash Z'
\\
@V{\a}VV @VV{\b}V
\\
\PP^1 @<{\cong}<< \PP^1.
\end{CD}
\end{equation}


\begin{thebibliography}{10}

\bibitem{AZ}
M. Artin and J.J. Zhang,
Non-commutative Projective Schemes,
{\it Adv. Math.}, {\bf 109} (1994), 228-287.

\bibitem{BGK}
V. Baranovsky, V. Ginzburg, and A. Kuznetzov, 
Quiver varieties and a noncommutative $\PP^2$,
Preprint (2001) AG/0103068.

\bibitem{ADR}
A. del Rio,
Graded rings and equivalences of categories,
{\it Comm. Alg.,} {\bf 19} (1991) 997-1012.

\bibitem{BR}
M. Beltrametti and L. Robbiano, Introduction to the theory
of weighted projective spaces, {\it Expo. Math.,} {\bf 4}
(1986) 111-162.


\bibitem{Dolg}
I. Dolgachev, Weighted projective varieties, in {\it Group
actions and vector fields,} pp. 34-71, Lecture Notes in Math. 173,
Springer-Verlag, 1973.


\bibitem{Gab}
P. Gabriel, Des Cat\'egories Ab\'eliennes, {\it Bull. Soc. Math. Fr.,}
{\bf 90} (1962) 323-448.

\bibitem{KKO}
A. Kapustin, A. Kuznetzov, and D. Orlov,
Noncommutative Instantons and Twistor Transform,
{\it Comm. Math. Phys.,} {\bf 221} (2001) 385--432


\bibitem{R}
A.L. Rosenberg, {\it Non-commutative algebraic geometry and
representations of quantized algebras,} Vol. 330, Kluwer
Academic Publishers, 1995.

\bibitem{Sm1}
S.P. Smith,
Subspaces of non-commutative spaces,
{\it Trans. Amer. Math. Soc.,} {\bf 354} (2002) 2131--2171.

\bibitem{Sm2}
S.P. Smith, Integral non-commutative spaces, {\it J.  Algebra,} 
{\bf 246} (2001) 793-810.



\bibitem{S}
D. Stephenson, Quantum planes of weight $(1,1,n)$,
{\it J. Algebra,} {\bf 225} (2000) 70-92.
Corrigendum, {\it J. Algebra,} {\bf 234} (2000) 277-278.

\bibitem{V}
A.B. Verevkin, On a non-commutative analogue of the
category of coherent sheaves on a projective scheme, {\em Amer. Math.
Soc.  Trans.,} {\bf 151} (1992) 41-53.

\bibitem{vdB}
M. Van den Bergh,
Blowing up of a non-commutative smooth surface, {\it Mem.
Amer.  Math. Soc.,} {\bf 154} (2001) no. 734.

 
\end{thebibliography}
\end{document}